\documentclass[12pt,american,twoside,a4wide]{article}
\usepackage[T1]{fontenc}
\usepackage[latin1]{inputenc}
\usepackage{latexsym}
\usepackage{amsmath}
\usepackage{bbm}
\usepackage{amssymb}
\usepackage{babel}
\usepackage{amsfonts}
\usepackage{times}
\usepackage{theorem}
\usepackage{epsfig}
\usepackage{enumerate}
\usepackage{color}
\usepackage[active]{srcltx}
\usepackage[colorlinks=true]{hyperref}
\usepackage{nameref}
\usepackage{fancyhdr}

\makeatletter
\let\orgdescriptionlabel\descriptionlabel
\renewcommand*{\descriptionlabel}[1]{%
  \let\orglabel\label
  \let\label\@gobble
  \phantomsection
  \edef\@currentlabel{#1}%
  \let\label\orglabel
  \orgdescriptionlabel{#1}%
}
\makeatother

\numberwithin{equation}{section}

\newcounter{smallarabics}

\newcounter{smallroman}

\newcommand{\ben}{\begin{enumerate}[{\rm (1)}]}
\newcommand{\een}{\end{enumerate}}

\newtheorem{theoreme}{Theorem}[section]
\newtheorem{proposition}[theoreme]{Proposition}
\newtheorem{lemma}[theoreme]{Lemma}
\newtheorem{definition}[theoreme]{Definition}

\def\rr{{\mathbb R}}
\def\zz{{\mathbb Z}}
\def\cc{{\mathbb C}}

\def\textsl{{}}

\def\Im{{\rm Im}\,}

\def\ch{{\frak h}}

\def\cT{\mathcal{T}}
\newcommand{\slim}{\mathop{\mathrm{s-lim}}\limits}

\def\c0inf{C_0^\infty}
\def\bep{\begin{proposition}}
\def\eep{\end{proposition}}

\def\proof{\noindent {\bf Proof.}\ \ }

\def\cR{{\cal R}}

 \def\cB{{\cal B}}

\def\f{{\rm f}}

\def\crys{{\rm crystal}}

\def\i{{\rm i}}
\newcommand{\beq}{\begin{equation}}
\newcommand{\eeq}{\end{equation}}
\newcommand{\bear}[1]{\begin{array}{#1}}
\newcommand{\ear}{\end{array}}

\def\sign{{\rm sign}}

\def\sp{{\hat e}}

\newcommand{\e}{\mathrm{e}}
\renewcommand{\i}{\mathrm{i}}

\renewcommand{\d}{\mathrm{d}}

\def\dint{\mathop{\int^\oplus}\limits}
\def\qed{$\Box$\medskip}

\def\cJ{{\cal J}}
\def\cT{{\cal T}}

\def\cI{{\cal I}}

\def\cO{{\cal O}}

\def\Th{{\rm Th}}
\def\bel{\begin{lemma}}
\def\eel{\end{lemma}}
\def\bet{\begin{theoreme}}
\def\eet{\end{theoreme}}
\def\bed{\begin{definition}}
\def\eed{\end{definition}}
\def\bar{\overline}

\def\12{\frac{1}{2}}

\def\e{{\rm e}}

\def\d{{\rm d}}
\def\Ran{{\rm Ran}\,}

\def\f{{\rm fr}}

\def\ac{{\rm ac}}

\def\sp{{\rm sp}}
\def\cS{{\cal S}}
\def\cR{{\cal R}}

\def\fh{{\mathfrak h}}

\def\tr{{\rm tr}}

\begin{document}
\def\today{}
\title{Landauer-B\"uttiker and Thouless conductance}
\author{L. Bruneau$^{1}$, V. Jak\v{s}i\'c$^{2}$, Y. Last$^3$, C.-A. Pillet$^4$
\\ \\ 
$^1$ D\'epartement de Math\'ematiques and UMR 8088\\
CNRS and Universit\'e de Cergy-Pontoise\\
95000 Cergy-Pontoise, France
\\ \\
$^2$Department of Mathematics and Statistics\\ 
McGill University\\
805 Sherbrooke Street West \\
Montreal,  QC,  H3A 2K6, Canada
\\ \\
$^3$Institute of Mathematics\\
The Hebrew University\\
91904 Jerusalem, Israel
\\ \\
$^4$Universit\'e de Toulon, CNRS, CPT, UMR 7332, 83957 La Garde, France\\
Aix-Marseille Universit\'e, CNRS, CPT, UMR 7332, 13288 Marseille, France}
\maketitle
\thispagestyle{empty}
\bigskip
\bigskip

\centerline {\large \bf  Dedicated to the memory of Markus B\"uttiker}

\bigskip
\newpage
{\small\noindent{\bf Abstract.} In the independent electron approximation, the 
average (energy/charge/entropy) current flowing through a finite sample $\cS$ 
connected to two electronic reservoirs can be computed by scattering theoretic 
arguments which lead to the famous Landauer-B\"uttiker formula. Another well 
known formula has been proposed by Thouless on the basis of a  scaling argument. The Thouless formula relates the 
conductance of the sample to the width of the 
spectral bands of the infinite crystal obtained by periodic juxtaposition of 
$\cS$. In this spirit, we define Landauer-B\"uttiker crystalline currents by extending  the Landauer-B\"uttiker formula to a setup 
where the  sample $\cS$ is replaced by  a  periodic structure whose unit cell is $\cS$.
We argue that these crystalline currents are closely related to the Thouless 
currents. For example, the crystalline heat current is bounded above by the 
Thouless heat current, and this bound saturates iff the coupling between the 
reservoirs and the sample is reflectionless. Our analysis leads to a rigorous 
derivation of the Thouless formula from the first principles of quantum 
statistical mechanics. 
}

\section{Introduction}
Recent rigorous formulations of the Landauer-B\"uttiker
formula~\cite{AJPP,BSP,CJM,N}
in the framework of nonequilibrium quantum statistical mechanics have opened the 
way to the mathematical study of a variety of related transport phenomena
in quantum mechanics~\cite{GJW,JLPa,JLPi}. This paper is a continuation of this 
line of research. Our main goal here is to provide a mathematically rigorous 
proof of the celebrated Thouless conductance formula.

This paper is organized as follows. In Section~\ref{sec-ebb} we review the 
Electronic Black Box Model and the corresponding Landauer-B\"uttiker formula. 
These topics have been discussed from both technical and pedagogical point of view 
in~\cite{AJPP,BSP} and the reader may consult these works for additional 
information. The notions of Thouless energy and Thouless conductance are reviewed in 
Section~\ref{sec-thouless}. Our main results are stated in Section~\ref{sec-CLB}. 
The proofs are given in Section~\ref{sec-proof}.

\bigskip\noindent
{\bf Acknowledgment.} The research of V.J. was partly supported by NSERC. The research of Y.L.\ was partly supported by The Israel Science Foundation
(Grant No.\ 1105/10) and by Grant No.\ 2010348 from the United States-Israel
Binational Science Foundation (BSF), Jerusalem, Israel. A part 
of this work has been done during a visit of L.B. to McGill University supported 
by NSERC. Another part was done during the visit of V.J. to The Hebrew University supported by ISF.

\subsection{The Electronic Black Box Model and  the Landauer-B\"uttiker formula}
\label{sec-ebb}
 
The Electronic Black Box Model (abbreviated EBBM) describes a finite sample $\cS$
connected to two infinitely extended electronic reservoirs $\cR_{l/r}$ (where 
$l/r$ stands for left/right) in the independent electron approximation. The 
coupling between the sample and the reservoirs allows for a flow of 
energy/charge/entropy through the joint system $\cR_l+\cS+\cR_r$. 

We shall restrict our attention to one-dimensional samples in the tight 
binding approximation. Thus, the sample $\cS$ is a free Fermi gas with 
one-particle Hilbert space $\ch_\cS=\ell^2(Z_L)$, where $Z_L=[1,L]\cap\zz$ is a 
finite lattice. Its Hamiltonian $h_\cS$ is a Jacobi matrix with parameters $\{J_x\}_{1\le x< L}$,
$\{\lambda_x\}_{1\le x\le L}$, 
\beq\label{jacobi-sample}
(h_\cS u)(x)=J_x u(x+1)+J_{x-1}u(x-1)+\lambda_x u(x),\qquad
x\in Z_L,
\eeq
and  Dirichlet boundary condition $u(0)=u(L+1)=0$.
The reservoir $\cR_{l/r}$ is a  free Fermi gas with one-particle Hilbert space
$\fh_{l/r}$ and one-particle Hamiltonian $h_{l/r}$. The one-particle Hilbert
space and Hamiltonian of the composite system $\cR_l+\cS+\cR_r$ are
\begin{align*}
\fh&=\fh_l\oplus\fh_\cS\oplus\fh_r,\\[3mm]
h_0&=h_l\oplus h_\cS \oplus h_r.
\end{align*}
The coupling of the sample with the reservoir $\cR_{l/r}$ is realized by the 
hopping Hamiltonian
\[
v_{l/r}=|\chi_{l/r}\rangle\langle\psi_{l/r}|+|\psi_{l/r}\rangle\langle\chi_{l/r}|,
\]
where $\chi_{l/r}\in\fh_{l/r}$ is a unit vector while $\psi_l=\delta_1$ 
and $\psi_r=\delta_L$ are Kronecker delta functions in $\ch_\cS$.
The one-particle Hamiltonian of the coupled system is 
\[
h=h_0+\kappa v = h_0 +\kappa ( v_l+v_r),
\]
where $\kappa\not=0$ is the coupling strength.
As observed in~\cite{BJP}, for the purposes of discussing transport properties of
the coupled system $\cR_l+\cS+\cR_r$ one may assume, without loss of generality,
that $\chi_{l/r}$ is a cyclic vector for $h_{l/r}$. Hence, passing to the spectral
representation we may assume that $h_{l/r}$ acts as multiplication by $E$ on
\[
\fh_{l/r}=L^2(\rr,\d\nu_{l/r}(E)),
\]
where $\nu_{l/r}$ is the spectral measure of $h_{l/r}$ associated to 
$\chi_{l/r}$. Moreover, in this representation one has $\chi_{l/r}(E)=1$
and $\nu_{l/r}(\rr)=1$.

Let $\Gamma_-(\fh)$ be the fermionic Fock space over $\fh$ and denote by
$a(f)/a^\ast(f)$ the annihilation/creation operator on $\Gamma_-(\fh)$ 
associated to $f\in\fh$. The algebra of observables of the EBBM is the 
$C^\ast$-algebra $\cO$ of all bounded operators on $\Gamma_-(\fh)$ generated by the identity $I$ and 
$\{a^\ast(f)a(g)\,|\,f,g\in\fh\}$.

Let $H=\d\Gamma(h)$ be the second quantized Hamiltonian. The map 
\[
\tau^t(a^\ast(f)a(g))=\e^{\i tH}a^\ast(f)a(g)\e^{-\i tH}
=a^\ast(\e^{\i th}f)a(\e^{\i th}g),
\]
extends to a strongly continuous group of $\ast$-automorphisms of $\cO$ and the 
pair $(\cO,\tau^t)$ is a $C^\ast$-dynamical system. The states of the EBBM are 
normalized positive linear functionals on $\cO$. Since we are dealing with 
independent electrons, quasi-free states on $\cO$ will be of particular
relevance. Let $\rho$ be a self-adjoint operator of $\fh$ such that
$0\leq\rho\leq I$. The gauge-invariant quasi-free state of density $\rho$ is the 
unique state $\omega_\rho$ on $\cO$ satisfying
\[
\omega_\rho(a^\ast(f_1)\cdots a^\ast(f_n)a(g_m)\cdots a(g_1))
=\delta_{n,m}\det\{\langle g_i,\rho f_j\rangle\},
\]
for any integers $n,m$ and all $f_1,\ldots,f_n,g_1,\ldots,g_m\in\fh$.

The initial state $\omega_0$ of the EBBM is the gauge-invariant quasi-free state of 
density
\[
\rho_0=\rho_l\oplus\rho_\cS\oplus\rho_r,
\]
where $\rho_{l/r}$ is the operator of multiplication by the Fermi-Dirac density 
\[
\rho_{l/r}(E)=\frac{1}{1+\e^{\beta_{l/r}(E-\mu_{l/r})}}.
\]
In other words, the reservoir $\cR_{l/r}$ is initially in thermal 
equilibrium at inverse temperature $\beta_{l/r}>0$ and chemical potential 
$\mu_{l/r}\in\rr$. None of our results depends on the particular choice of
the initial density $\rho_\cS$ of the sample.

The EBBM is the quantum dynamical system $(\cO,\tau^t,\omega_0)$.

The basic questions regarding this system concern its behavior in the large time
limit $t\to\infty$. To deal with this limit we observe that the initial state 
$\omega_0=\omega_{\rho_0}$ evolves in time as
\[
\omega_t=\omega_0\circ\tau^t=\omega_{\rho_t},
\]
where the density at time $t$ is given by
\[
\rho_t=\e^{-\i th}\rho_0\e^{\i th}.
\]
One expects that $\omega_t\to\omega_+$ as $t\rightarrow \infty$, where the 
non-equilibrium steady state $\omega_+$ carries energy/charge/entropy current
induced by the initial temperature/chemical potential differential. These 
currents can be expressed by Landauer-B\"uttiker formulas involving
transmission properties of individual electrons evolving under the dynamics 
generated by the one-particle Hamiltonian $h$. These heuristics is mathematically 
formalized as follows. 

The observables describing energy and charge current out of $\cR_{l/r}$ are 
\begin{align*}
\Phi_{l/r}&=\d\Gamma(-\i[h,h_{l/r}])
=\kappa\left(a^\ast(\i h_{l/r}\chi_{l/r})a(\psi_{l/r})
+a^\ast(\psi_{l/r})a(\i h_{l/r}\chi_{l/r})\right),\\[3mm]
\cI_{l/r}&=\d\Gamma(-\i[h,1_{l/r}])
=\kappa\left(a^\ast(\i\chi_{l/r})a(\psi_{l/r})
+a^\ast(\psi_{l/r})a(\i\chi_{l/r})\right).
\end{align*}
The entropy current associated to the heat flux dissipated into the reservoirs is
\[
\cJ=-\beta_l(\Phi_l-\mu_l\cI_l)-\beta_r(\Phi_r-\mu_r\cI_r).
\]
The transmittance of the sample is the function \footnote{This function is defined for Lebesgue a.e. $E\in \rr$.}
\beq\label{tra-coef}
\rr\ni E\mapsto
\cT(E)=4\kappa^4|\langle\psi_l,(h-E-\i0)^{-1}\psi_r\rangle|^2\,\Im F_l(E)\,\Im F_r(E),
\eeq
where  
\beq\label{reserv-funct}
F_{l/r}(E)=\langle\chi_{l/r},(h_{l/r}-E-\i0)^{-1}\chi_{l/r}\rangle.
\eeq
Note that $\Im F_{l/r}(E)\ge0$ and \footnote{This identity is 
understood modulo sets of Lebesgue measure zero.} 
\beq\label{non-tri}
\{E\in\rr\,|\,\cT(E)>0\}=\Sigma_{l}\cap \Sigma_{r},
\eeq
where 
\beq
\Sigma_{l/r}=\{E\in\rr\,|\,\Im F_{l/r}(E)>0\}
\label{ImF}
\eeq
is the essential support of the absolutely continuous spectrum of $h_{l/r}$.

For $E\in\rr$ we set
$$
\zeta_{l/r}(E)=\beta_{l/r}(E-\mu_{l/r}),\qquad
\Delta_{l/r}(E)=\rho_{l/r}(E)-\rho_{r/l}(E),
$$
and
$$
\varsigma(E)=(\zeta_{r}(E)-\zeta_{l}(E))\Delta_l(E)
=(\zeta_{l}(E)-\zeta_{r}(E))\Delta_r(E).
$$
Note that if $\beta_l=\beta_r$ and $\mu_l=\mu_r$ then $\varsigma(E)$ 
vanishes identically but that it is strictly positive for all $E\in\rr$ otherwise.

Our basic assumption in the following is:
\begin{quote}{\bf Assumption A.} The one-particle Hamiltonian $h$ has no 
singular continuous spectrum.
\end{quote}
The starting point of this paper is the following result~\cite{AJPP}.
\bet\label{basic}
Suppose that Assumption ${\bf A}$ holds. Then for all $A\in\cO$ the limit 
\[
\omega_+(A)=\lim_{t\to\infty}\frac1t\int_0^t\omega_s(A)\d s,
\]
exists and defines a state $\omega_+$ on $\cO$. Moreover, the following 
formulas for the average steady currents hold:
\beq
\begin{split}
\langle\Phi_{l/r}\rangle_+&=\omega_+(\Phi_{l/r})
=\frac{1}{2\pi}\int_\rr\cT(E)E\Delta_{l/r}(E)\d E,\\[3mm]
\langle\cI_{{l/r}}\rangle_+&=\omega_+(\cI_{l/r})
=\frac{1}{2\pi}\int_\rr\cT(E)\Delta_{l/r}(E)\d E,\\[3mm]
\langle\cJ\rangle_+&=\omega_+(\cJ)
=\frac{1}{2\pi}\int_\rr\cT(E)\varsigma(E)\d E.
\end{split}
\label{landauer-buttiker}
\eeq

\eet
We finish this section with a number of remarks regarding this result. 

\bigskip
{\bf\noindent Remark 1.} Theorem~\ref{basic} deals with the simplest non-trivial
setting in the study of electronic transport in the independent electron approximation. 
Various generalizations of this setting and of Theorem~\ref{basic} can be found 
in~\cite{AJPP,N,BSP}.

{\bf\noindent Remark 2.} We have chosen units in such a way that the electronic 
charge $e$, the reduced Planck constant $\hbar$, and the Boltzmann constant $k_B$ are
unity. The energy, charge, and entropy currents in the above formulas are 
expressed in units of $1/\hbar$, $e/\hbar$ and $k_B/\hbar$. Note also that these 
formulas do not include the spin degeneracy which should be accounted for
by a factor $2$.

{\bf\noindent Remark 3.} If Assumption {\bf A} is replaced by the stronger 
assumption that $h$ has no singular spectrum, then 
\[
\omega_+(A)=\lim_{t\to\infty}\omega_t(A),
\]
holds for all $A\in\cO$. 

{\bf\noindent Remark 4.} If $\beta_l=\beta_r=\beta$ and $\mu_l=\mu_r=\mu$, then $\omega_+$
is the thermal equilibrium state of the coupled system $\cR_l+\cS+\cR_r$ for the
given intensive thermodynamic parameters, i.e., the quasi-free state of density
$$
\rho_{\beta,\mu}=\frac1{1+\e^{\beta(h-\mu)}}.
$$
In this case, all currents vanish in average. In the following we shall exclude
this possibility and assume that $\beta_l\not=\beta_r$ or/and
$\mu_l\not=\mu_r$. The state $\omega_+$ is then a non-equilibrium steady state 
of the EBBM and~\eqref{landauer-buttiker} are the Landauer-B\"uttiker
formulas for the steady state  currents.

{\bf\noindent Remark 5.} Apart from the choice of the intensive thermodynamic 
parameters $\beta_{l/r}$ and $\mu_{l/r}$, the transmission coefficient $\cT(E)$ 
completely determines the steady state currents. It follows from~\eqref{non-tri} 
that the steady state currents are non-vanishing iff the Lebesgue measure  $|\Sigma_l\cap\Sigma_r|$ of the 
set $\Sigma_l\cap\Sigma_r$ is strictly positive, i.e., iff there exists an open 
scattering channel between the left and the right reservoir. Note the obvious 
conservation laws
$$
\omega_+(\Phi_l)+\omega_+(\Phi_r)=0,\qquad
\omega_+(\cI_l)+\omega_+(\cI_r)=0.
$$

Entropy balance in the steady state $\omega_+$ implies that
\[
\langle\cJ\rangle_+
=-\beta_l(\langle\Phi_l\rangle_+-\mu_l\langle\cI_l\rangle_+)
 -\beta_r(\langle\Phi_r\rangle_+-\mu_r\langle\cI_r\rangle_+),
\]
coincides with the rate of entropy production in $\cS$ \cite{AJPP}. It follows from the
above observation that $\langle\cJ\rangle_+\ge0$, and that the  inequality is strict
iff $|\Sigma_l\cap\Sigma_r|>0$.

{\bf\noindent Remark 6.} The proof of Theorem~\ref{basic} is based on the
scattering theory of the pair $(h,h_0)$ and elucidates the physical meaning  
of $\cT(E)$. It follows from the trace class scattering theory that the wave 
operators \footnote{$1_\ac(h_0)$ is the spectral projection on the absolutely 
continuous part of the spectrum of $h_0$.}
\[
w_{\pm}=\slim_{t\to\pm \infty}\e^{\i th}\e^{-\i th_0}1_{\ac}(h_0),
\]
exist and are complete. The scattering matrix $s=w_+^\ast w_-$ is a unitary 
operator on 
\[\fh_\ac(h_0)=\Ran 1_{\ac}(h_0)=\Ran 1_\ac(h_l)\oplus\Ran 1_\ac(h_r),\]
 and acts as the operator of multiplication by a unitary $2\times 2$ matrix 
\[
s(E)=\left[\begin{matrix}
s_{ll}(E)&s_{lr}(E)\\
s_{rl}(E)&s_{rr}(E)
\end{matrix}\right].
\]
One then has 
\[
\cT(E)=|s_{lr}(E)|^2=|s_{rl}(E)|^2,
\]
i.e., the transmittance $\cT(E)$ is the transmission probability between the left and 
right reservoir at energy $E$.

\medskip
{\bf\noindent Remark 7.} An additional insight into $\cT(E)$ can be obtained 
by taking into account the spatial structure of the reservoirs. Suppose that the 
reservoirs are non-trivial in the sense that the measures $\nu_{l/r}$ have 
non-vanishing absolutely continuous component. The standard orthogonal polynomial 
construction (see Theorem I.2.4 in \cite{Si}) provides a unitary operator
$U:\fh\to\ell^2(\zz)$ such that the following holds:
\begin{enumerate}[{\sl (i)}]
\item $U\fh_l=\ell^2(]-\infty,0]\cap\zz)$, $U\fh_\cS=\fh_\cS$ and
$U\fh_r=\ell^2([L+1,\infty[\cap\zz)$.
\item There is a Jacobi matrix on $\ell^2(]-\infty,0]\cap\zz)$, with parameters 
$\{J_x\}_{x<0}$, $\{\lambda_x\}_{x\le0}$ and Dirichlet boundary condition such that
for $u\in\ell^2(]-\infty,0]\cap\zz)$
$$
(Uh_{l}U^\ast u)(x)=J_xu(x+1)+J_{x-1}u(x-1)+\lambda_{x}u(x),\qquad (x\le0, u(1)=0).
$$
\item $Uh_\cS U^\ast=h_\cS$.
\item  There is a Jacobi matrix on $\ell^2([L+1,\infty[\cap\zz)$, with parameters 
$\{J_x\}_{x> L}$, $\{\lambda_x\}_{x>L}$ and Dirichlet boundary condition such that
for $u\in\ell^2([L+1,\infty[\cap\zz)$
$$
(Uh_{r}U^\ast u)(x)=J_xu(x+1)+J_{x-1}u(x-1)+\lambda_{x}u(x),\qquad (x>L, u(L)=0).
$$
\item $U\chi_l=\delta_{0}$ and $U\chi_r=\delta_{L+1}$.
\item If $J_0=J_L=\kappa$,  then for $u \in \ell^2(\zz)$
$$
(UhU^\ast u)(x)=J_{x}u(x+1)+J_{x-1}u(x-1)+\lambda_x u(x).
$$
\end{enumerate}

It follows that our  EBBM is unitarily equivalent to a Jacobi matrix EBBM on $\ell^2(\zz)$.

For $z\in\cc_+$ let $u_{l/r}(\,\cdot\,,z)$ be the unique solution of the equation 
\beq
J_x u_{l/r}(x+1,z )+J_{x-1}u_{l/r}(x-1,z)+\lambda_x u_{l/r}(x,z)=zu_{l/r}(x,z),
\label{well}
\eeq
that is square summable at $\pm\infty$ and normalized by $u_{l/r}(0,z)=1$.
For all $x\in\zz$ and Lebesgue a.e.\;$E\in \rr$ the limit 
\[
\lim_{\epsilon\downarrow0}u_{l/r}(x,E+\i\epsilon)=u_{l/r}(x,E),
\]
exists, is finite, and solves~\eqref{well} with $z=E$. For Lebesgue 
a.e.\;$E\in\Sigma_r$ the solution $u_r(E)$ is not a multiple of a real 
solution and so $\bar{u_{r}(E)}$ is also a solution of~\eqref{well} 
linearly independent of $u_r(E)$. Hence, for Lebesgue a.e.\;$E\in\Sigma_r$ 
one has
\[
u_l(E)=\alpha (E)\bar{u_r(E)}+\beta(E)u_r(E).
\]
The spectral reflection probability of~\cite{GNP,GS} is 
\[
R_r(E)=\left|\frac{\beta(E)}{\alpha(E)}\right|^2.
\]
One extends $R_r$ to $\rr$ by setting $R_r(E)=1$ for $E\not\in\Sigma_r$. 
$R_l(E)$ is defined analogously. A simple computation (see Section~5 
in~\cite{JLPa}) gives
\[
R_r(E)=R_l(E)=|s_{ll}(E)|^2=|s_{rr}(E)|^2.
\]
Hence, 
\[
\cT(E)=1-R_r(E)=1-R_l(E),
\]
is the spectral transmission probability at the energy $E$.

\subsection{The Thouless energy and the Thouless conductance formula}
\label{sec-thouless}

We start with a review of the notions of Thouless energy and conductance 
as discussed in the physics literature. We follow~\cite{La} and our main goal is
to extract mathematically well-defined quantities that correspond to these 
heuristic notions.

\subsubsection{Heuristics}

Let $\delta t$ be the typical time spent in $\cS$ by an 
electron on its journey from one reservoir to the other. The Mandelstam-Tamm 
time-energy uncertainty relation $\delta E\,\delta t\gtrsim1$ (\cite{MT}, see also \cite{FP,Bu}) 
 imposes a lower bound on the energy spread $\delta E$
of the electron wave function. This sets the {\sl Thouless energy\/} scale 
$E_\Th\sim\delta E$ \cite{Th}. Assuming a diffusive behavior, one has
$L^2\sim D\delta t$ where $L$ is the sample size and $D$ the diffusion constant, 
and hence
$$
E_\Th\gtrsim\frac D{L^2}.
$$
Einstein's relation $$
\sigma=D\varrho\lesssim L^2E_\Th\varrho,
$$
further links $D$ to the conductivity $\sigma$ of the sample
and its density of states $\varrho$.
Finally, $\varrho$ relates to the typical energy level spacing $\Delta E$ of the 
sample as 
$$
\varrho L\Delta E\sim 1.
$$
It follows that
$$
\sigma\lesssim L\frac{E_\Th}{\Delta E},
$$
and hence that the conductance $g=L^{-1}\sigma$ of the one-dimensional sample 
satisfies
\beq
g\lesssim g_\Th=\frac{E_\Th}{\Delta E}.
\label{ThBound}
\eeq
The quantity at the r.h.s.\;of this formula is known as the {\sl Thouless 
conductance.} Although the above derivation is  extremely 
heuristic in nature, the Thouless conductance $g_\Th$ and the closely related 
Thouless energy $E_\Th$ are  widely accepted by both 
theoretical and experimental physicists and play an important role in
the scaling theory of localization \cite{AALR}. 

The above argument and the resulting inequality on the l.h.s.\;of \eqref{ThBound}
suggest that we should consider $g_{\rm Th}$
as an upper bound on the conductance of the sample $\cS$. Indeed, for such a 
microscopic system the conductance is not an intrinsic property of the sample, but 
depends also on  the reservoirs and the nature of the coupling, which 
determines the  reservoir's ability to feed the available energy levels of the sample. Thus, one
expects that  saturation of the Thouless bound \eqref{ThBound} occurs for optimal feeding of 
the sample by the reservoirs, a property of the joint system $\cR_l+\cS+\cR_r$ which we 
shall try to elucidate in the remaining part of this paper (see 
Remarks 1 and 4 after Theorem~\ref{final-thouless}).

\subsubsection{The crystal model}
\label{Sect:crystallineEBBM}

There is a simple way to make the sample transparent to incoming electrons and hence
to ensure its optimal feeding: it suffices to implement the reservoirs in such a way 
that the joint system $\cR_l+\cS+\cR_r$ is periodic. We shall call {\sl crystalline
EBBM\/} the model obtained by repeating the sample so as to obtain a periodic 
crystal with unit cell $\cS$ (a construction closely related to the scaling 
argument of \cite{ET}).

Consider the periodic Jacobi matrix $h_\crys$ on
$\ell^2(\zz)$ obtained by extending the Jacobi parameters $\{\lambda_x\}_{1\le x\le L}$
and $\{J_x\}_{1\le x<L}$ of the sample Hamiltonian
$h_\cS$ to the entire lattice $\zz$ by setting $J_L=\kappa_\cS$ and 
$$
J_{x+nL}=J_x, \qquad \lambda_{x+nL}=\lambda_x,
$$
for any $n\in\zz$ and $x\in Z_L$. The internal coupling constant $\kappa_\cS$ is a priori
an arbitrary parameter, except for the obvious constraint $\kappa_\cS\not=0$. In 
practice it will be determined by the physics of the problem. In a model where the 
sample parameters $J_x$ are independent copies of a random variable, $\kappa_\cS$ will 
be another instance of this variable. If $h_\cS$ is a discrete Schrödinger operator
with $J_x=J$ for all $x\in Z_L$ the choice of $\kappa_\cS=J$ appears natural.

In the crystaline  EBBM model, the single particle Hilbert spaces of 
the reservoirs are $\fh_l=\ell^2(]-\infty,0]\cap\zz)$ and 
$\fh_r=\ell^2([L+1,\infty[\cap\zz)$, the corresponding single particle Hamiltonians are Jacobi matrices  with parameters 
 $\{\{J_x\}_{x<0}, \{\lambda_x\}_{x<0}\}$, $\{\{J_x\}_{x>L}, \{\lambda_x\}_{x>L}\}$ and Dirichlet boundary condition, 
 $\chi_l=\delta_0$, $\chi_r=\delta_{L+1}$, and the coupling constant is set to $\kappa=\kappa_\cS$. 
 The one particle Hamiltonian of the coupled system is $h_\crys$.
 
 We emphasize that the crystallization of the sample is specified by the pair $(\cS, \kappa_\cS)$ and not by $\cS$ alone. 

\begin{figure}
\centering
\includegraphics[scale=1]{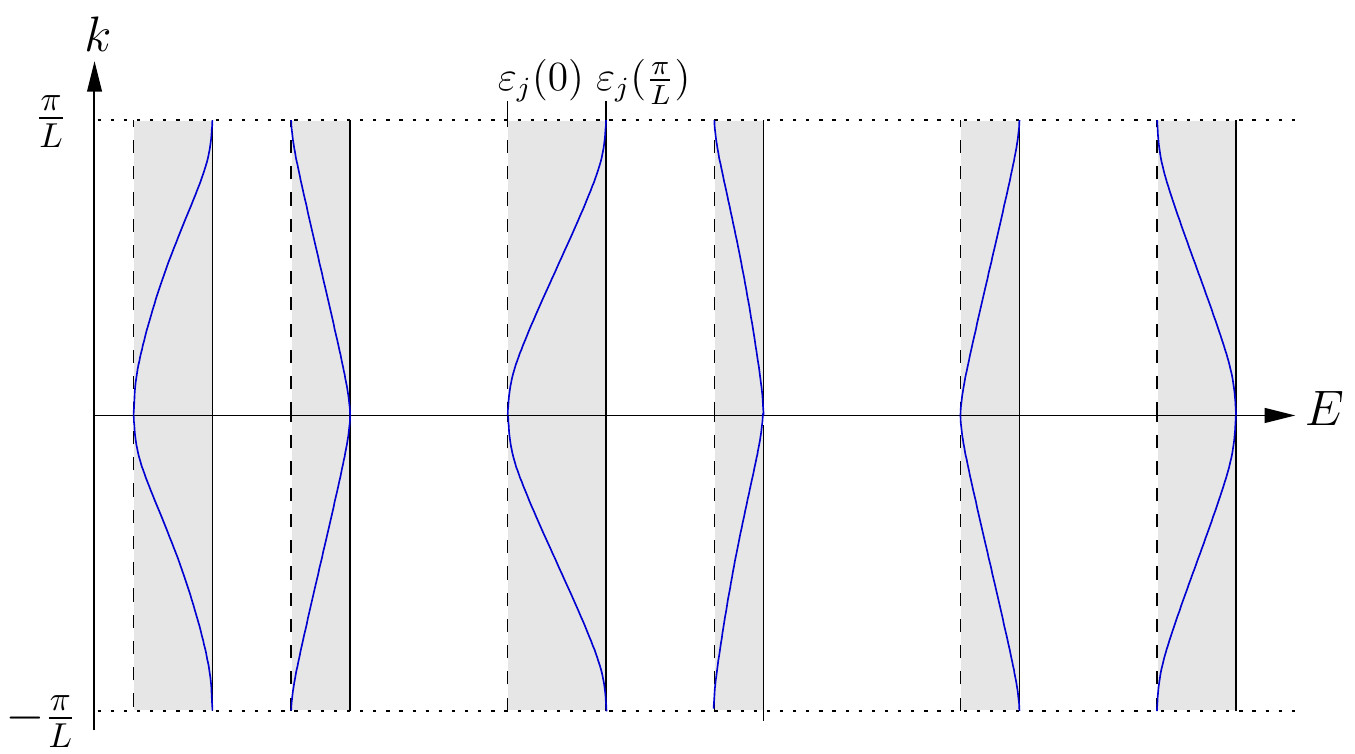}
\caption{Eigenvalues of $h(k)$ and the band spectrum of $h_\crys$.}
\label{Fig1}
\end{figure}

The Bloch-Floquet decomposition of $h_\crys$ reads
$$
\ell^2(\zz)=\dint_{\mathfrak{B}_L}\fh_\cS\,\d k,\qquad
h_\crys=\dint_{\mathfrak{B}_L}h(k)\,\d k,
$$
where $\mathfrak{B}_L=[-\pi/L,\pi/L]$ is the first Brillouin zone of the crystal and
$h(k)$ is the self-adjoint operator on $\fh_\cS$ obtained from~\eqref{jacobi-sample} 
by replacing Dirichlet by Bloch boundary conditions
$$
u(0)=\e^{-\i k L}u(L),\qquad u(L+1)=\e^{\i k L}u(1).
$$
The spectrum of $h(k)$ consists of $L$ eigenvalues 
$\varepsilon_1(k)\le\cdots\le\varepsilon_L(k)$ which are even functions of $k$, 
real analytic and strictly monotone on $]0,\pi/L[$. Moreover
$$
\varepsilon_L(0)>\varepsilon_L(\pi/L)\ge\varepsilon_{L-1}(\pi/L)
>\varepsilon_{L-1}(0)\ge\varepsilon_{L-2}(0)>\cdots
$$
Thus, the spectrum of $h_\crys$ is
$$
\sp(h_\crys)=\bigcup_{k\in\cB_L}\sp(h(k))=\bigcup_{j=1}^L B_j,
$$
where $B_j$ is a closed interval with boundary points $\varepsilon_j(0)$ and $\varepsilon_j(\pi/L)$ (see Theorem 5.3.4 in~\cite{Si} and Figure~\ref{Fig1}).

For simplicity, suppose that the reservoirs are at zero temperature, i.e., 
$\beta_l=\beta_r=\infty$, and assume that $\mu_l<\mu_r$. In this case the
function $\Delta_r$ is the characteristic function of the interval $[\mu_l,\mu_r]$
and the Landauer-B\"uttiker formula for the steady charge current out 
of the right reservoir becomes
\[
\omega_+(\cJ_r)=\frac{1}{2\pi}\int_{\mu_l}^{\mu_r}\cT(E)\d E.
\]
Since the transmittance of a unit cell in a perfect crystal is the characteristic
function of its spectrum, we get
\[
\omega_+(\cJ_r)=\frac1{2\pi}|\sp(h_\crys)\cap [\mu_l,\mu_r]|,
\]
from which we infer that the mean conductance of the sample on the
energy window $I=[\mu_l,\mu_r]$ is given by
$$
g(I)=\frac1{2\pi}\frac{|\sp(h_\crys)\cap I|}{|I|}.
$$

\begin{figure}
\centering
\includegraphics[scale=0.5]{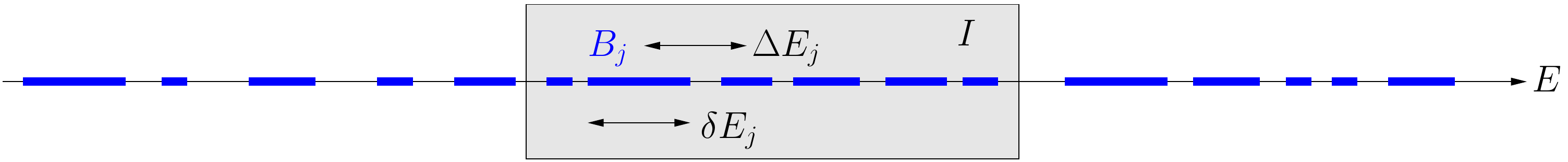}
\caption{Energy uncertainty $\delta E_j$ and level spacing $\Delta E_j$ for the $j$-th 
band in the window $I$.}
\label{Fig2}
\end{figure}

We shall now argue that this expression can be interpreted as the
Thouless conductance associated to the energy window $I$.
In order for the level spacing $\Delta E$ to be well defined, the interval $I$
should contain several bands $B_j$ of $h_\crys$. The energy uncertainty within a 
single band $B_j\subset I$ is of the order of the band width $|B_j|=|\varepsilon_j(\pi/L)-\varepsilon_j(0)|$ which coincides with the
variation of the eigenvalue $\varepsilon_j(k)$ as the Bloch boundary condition changes
from periodic to anti-periodic (see Figure~\ref{Fig2}). A rough but convenient 
estimate of the energy uncertainty within the window $I$ is given by the arithmetic 
mean
$$
\delta E\sim\frac{\sum_{B_j\subset I}|B_j|}{\sum_{B_j\subset I}1}
\sim\frac{|\sp(h_\crys)\cap I|}{\sum_{B_j\subset I}1}.
$$
On the other hand, the mean level spacing within $I$ is given by
$$
\Delta E\sim\frac{|I|}{\sum_{B_j\subset I}1},
$$
and the Thouless conductance becomes
$$
g_\Th\sim\frac{\delta E}{\Delta E}\sim\frac{|\sp(h_\crys)\cap I|}{|I|}.
$$

According to the previous arguments, we shall {\sl define} the Thouless
conductance of the pair $(\cS, \kappa_\cS)$ for the energy window $I\subset\rr$ by
\beq
g_\Th(I)=\frac1{2\pi}\frac{|\sp(h_\crys)\cap I|}{|I|},
\label{gTh}
\eeq
where $h_\crys$ is the periodization of $h_\cS$ with $J_L=\kappa_\cS$. 

Note that if $\kappa_\cS=0$, then $\sp(h_\crys)=\sp(h_\cS)$ and  $g_\Th(I)=0$ for all 
intervals $I$. More generally, 
the
width of the $L$ Bloch bands satisfies $|B_j|=\mathcal{O}(\kappa_\cS)$ as $\kappa_\cS\to0$, and 
$$
g_\Th(I)=\mathcal{O}(\kappa_\cS)
$$
in this limit.

\section{The crystalline limit}
\label{sec-CLB}

To further elaborate on the connection between Thouless conductance and
the Landauer-B\"uttiker formula, we shall now consider the approximation
of $h_\crys$ by finite repetitions of the sample $\cS$ connected to
arbitrary reservoirs.

\begin{figure}
\centering
\includegraphics[scale=0.4]{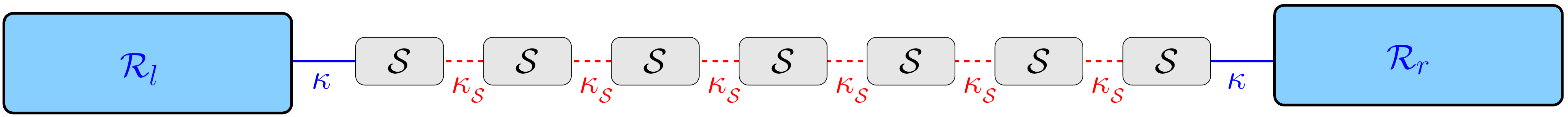}
\caption{The EBBM described by the Hamiltonian $h^{(N)}$ for $N=7$.}
\label{Fig3}
\end{figure}

Let $h_\crys$ be as in the previous section. Given a positive integer $N$, let 
$h_\cS^{(N)}$ be the restriction of $h_\crys$ to the finite lattice 
$Z_{NL}=[1,NL]\cap\zz$ with Dirichlet boundary condition. Hence $h_\cS^{(N)}$ 
is a Jacobi matrix acting on $\fh_\cS^{(N)}=\ell^2(Z_{NL})$ whose Jacobi 
parameters satisfy
$$
J_{x+nL}=J_x,\qquad \lambda_{x+nL}=\lambda_x,
\qquad (x\in Z_L,n=0,1,\ldots,N-1),
$$
where $\{J_x\}_{1\leq x<L}$ and $\{\lambda_x\}_{1\leq x\leq L}$ are the Jacobi 
parameters of the original sample Hamiltonian $h_\cS$ and $J_L=\kappa_\cS$. The pairs 
$(\fh_\cS^{(N)},h_\cS^{(N)})$ define a sequence of sample systems which are 
coupled to the reservoirs $\cR_{l/r}$ as in Section~\ref{sec-ebb}. The reservoirs' single particle Hilbert spaces and Hamiltonains
$(\ch_{l/r},h_{l/r})$, the vectors $\chi_{l/r}$, and the coupling strength $\kappa$ do not depend on $N$, and one takes 
$\psi_l=\delta_1$, $\psi_r=\delta_{NL}$. 
We assume that the one-particle 
Hamiltonian $h^{(N)}$ of the coupled systems satisfies Assumption~{\bf A} for all $N$ \footnote{Besides the crystaline ones, see Section \ref{Sect:crystallineEBBM}, a concrete example of reservoirs  where this is the case is $\ch_{l/r}=\ell^2(\zz_+)$, $h_{l/r}=-k\Delta$, $k>0$. For other examples and general results regarding this point we refer the reader to \cite{GJW}.}
and we denote by $\omega_+^{(N)}$, $\Phi_{l/r}^{(N)}$, $\cI_{l/r}^{(N)}$, 
$\cJ^{(N)}$ the respective NESS and flux observables. We are interested in 
the large $N$ limit of the charge, energy and entropy steady currents
(see Figure \ref{Fig3}).

Let $h_\crys^{(l)}$ and $h_\crys^{(r)}$ be the restrictions of $h_\crys$ to 
$\ell^2((-\infty,0]\cap\zz)$  and $\ell^{2}([1,\infty)\cap\zz)$ with Dirichlet
boundary conditions.  Denote by 
\beq
\begin{split}
m_l(E)&=\langle\delta_{0},(h_\crys^{(l)}-E-\i0)^{-1}\delta_{0}\rangle,\\[3mm]
m_r(E)&=\langle\delta_{1},(h_\crys^{(r)}-E-\i0)^{-1}\delta_{1}\rangle,
\end{split}
\label{Mfunct}
\eeq
the respective Weyl $m$-functions. One easily shows that $\Im m_{l/r}(E)>0$ for Lebesgue a.e. $E\in \sp(h_\crys)$. We set 
$\cT_\infty(E)=0$ for $E\in\rr\setminus(\sp(h_\crys)\cap\Sigma_l\cap\Sigma_r)$ and 
\beq
\cT_\infty(E)=\left[1+\frac14\left(
\frac{|\kappa_\cS^2m_r(E)-\kappa^2F_r(E)|^2}
{\Im(\kappa_\cS^2m_r(E))\Im(\kappa^2F_r(E))}
+\frac{|\kappa_\cS^2m_l(E)-\kappa^2F_l(E)|^2}
{\Im(\kappa_\cS^2m_l(E))\Im(\kappa^2F_l(E))}
\right)\right]^{-1}
\label{def:crys-tra-coef}
\eeq
for $E\in\sp(h_\crys)\cap\Sigma_l\cap\Sigma_r$. Obviously, $0\le\cT_\infty(E)\le1$ for Lebesgue a.e. $E\in \rr$.

Let 
\[\cT_N(E)=4\kappa^4|\langle\psi_l,(h^{(N)}-E-\i0)^{-1}\psi_r\rangle|^2\,\Im F_l(E)\,\Im F_r(E),
\]
be the transmittance of the $N$-fold 
repeated pair $(\cS, \kappa_\cS)$.
Our main technical result is:
\bet\label{prop:crys-tra-coef-bound}
For any $f\in L^1(\rr)$ one has
\beq
\lim_{N\to\infty}\int\cT_N(E)f(E)\,\d E=\int\cT_\infty(E)f(E)\,\d E.
\label{distribconverge}
\eeq

\eet
The proof of this theorem  is given in Section~\ref{sec-proof}.
As an immediate consequence, one has
\bet\label{crys-buttiker}
\beq\label{landauer-buttiker-crys}
\begin{split}
\langle\Phi_{l/r}\rangle_\infty
&=\lim_{N\to\infty}\omega_+^{(N)}(\Phi_{l/r}^{(N)})
=\frac{1}{2\pi}\int_{\sp(h_\crys)}\cT_\infty(E)E\Delta_{l/r}(E)\,\d E,\\[3mm]
\langle\cI_{l/r}\rangle_\infty
&=\lim_{N\to\infty}\omega_+^{(N)}(\cI_{l/r}^{(N)})
=\frac{1}{2\pi}\int_{\sp(h_\crys)}\cT_\infty(E)\Delta_{l/r}(E)\,\d E,\\[3mm]
\langle\cJ\rangle_\infty
&=\lim_{N\to\infty}\omega_+^{(N)}(\cJ^{(N)})
=\frac{1}{2\pi}\int_{\sp(h_\crys)}\cT_\infty(E)\varsigma(E)\,\d E.
\end{split}
\eeq
\eet
Obviously, the conservation laws
\begin{align*}
\langle\Phi_l\rangle_\infty+\langle\Phi_r\rangle_\infty&=0,\\[3mm]
\langle\cI_l\rangle_\infty+\langle\cI_r\rangle_\infty&=0,
 \end{align*}
hold, as well as the entropy balance relation 
\[
\langle\cJ\rangle_\infty
=-\beta_l(\langle\Phi_l\rangle_\infty-\mu_l\langle\cI_l\rangle_\infty)
 -\beta_r(\langle\Phi_r\rangle_\infty-\mu_r\langle\cI_r\rangle_\infty).
\]




\subsection{Thouless conductance revisited}
\label{sec-thouless-new}

Given a finite sample described by $\fh_\cS$, $h_\cS$, $\kappa_\cS$  and its periodization 
$h_\crys$, the Thouless currents associated to chemical potentials $\mu_{l/r}$ and 
inverse temperatures $\beta_{l/r}$ are defined by setting $\cT_\infty(E)=1$ in the 
formulas~\eqref{landauer-buttiker-crys}, i.e., by
assuming that the  transport between the reservoirs is reflectionless. The Thouless 
current formulas are:   
\beq\label{thoulesscurrents}
\begin{split}
\langle\Phi_{l/r}\rangle_\Th
&=\frac{1}{2\pi}\int_{\sp(h_\crys)}E\Delta_{l/r}(E)\d E,\\[3mm]
\langle\cI_{l/r}\rangle_\Th
&=\frac{1}{2\pi}\int_{\sp(h_\crys)}\Delta_{l/r}(E)\d E,\\[3mm]
\langle\cJ\rangle_\Th
&=\frac{1}{2\pi}\int_{\sp(h_\crys)}\varsigma(E)\d E.
\end{split}
\eeq
One has again the conservation laws
\begin{align*}
\langle\Phi_l\rangle_\Th+\langle\Phi_r\rangle_\Th&=0,\\[3mm]
\langle\cI_l\rangle_\Th+\langle\cI_r\rangle_\Th&=0,
 \end{align*}
and 
\[
\langle\cJ\rangle_\Th
=-\beta_l(\langle\Phi_l\rangle_\Th-\mu_l\langle\cI_l\rangle_\Th)
-\beta_r(\langle\Phi_r\rangle_\Th-\mu_r\langle\cI_r\rangle_\Th).
\]
Recall the definition of crystalline EBBM given in Section~\ref{Sect:crystallineEBBM}.
The following theorem is a direct consequence of the definition of Thouless 
currents~\eqref{thoulesscurrents} and of Theorem~\ref{crys-buttiker}.
\bet\label{final-thouless}
\ben 
\item 
\[
\langle\cJ\rangle_\Th=\sup\,\,\langle\cJ\rangle_\infty,
\]
where the supremum is taken over all realizations of the reservoirs. Moreover, this 
supremum is achieved if the EBBM is crystalline.
\item If $\mu_l=\mu_r=\mu$, $\mu\leq\inf\sp(h_\crys)$, and $\beta_l>\beta_r$, 
then 
\[
\langle\Phi_r\rangle_\Th=\sup\,\,\langle\Phi_r\rangle_\infty,
\]
and the supremum is achieved if the EBBM is crystalline. If $\beta_l<\beta_r$, 
then this result holds for $\langle\Phi_l\rangle$.  
If $\mu\geq\sup\sp(h_\crys)$, the same results hold with exchange of $l$ and $r$. 
\item If $\beta_l=\beta_r$ and $\mu_l<\mu_r$, then 
\[
\langle\cI_r\rangle_\Th=\sup\,\,\langle\cI_r\rangle_\infty,
\]
and the supremum is achieved if the EBBM is crystalline. If $\mu_l>\mu_r$, then 
this result holds for $\langle\cI_l\rangle$. 
\een
\eet

{\bf Remark 1.} All suprema in the previous theorem are achieved iff the transport 
between the reservoirs is reflectionless, that is, iff $\cT_\infty(E)=1$ for Lebesgue a.e. $E\in\sp(h_\crys)$. 
The crystalline EBBM provide such reservoirs. To elucidate this point further, note 
that $\cT_\infty(E)=1$ for Lebesgue a.e. $E\in\sp(h_\crys)$ iff 
$\sigma(h_\crys)\subset \Sigma_l\cap \Sigma_r$ and 
$$
\kappa_\cS^2m_{r/l}(E)=\kappa^2 F_{r/l}(E),
$$
for Lebesgue a.e. $E\in \sp(h_\crys)$. Since $\sp(h_\crys)$ has positive Lebesgue measure, 
the theory of boundary values of analytic functions (see \cite{Ja} or any book on harmonic analysis) yields that for 
any $z=E+\i \epsilon$ with $\epsilon >0$, 
\[\kappa_\cS^2\langle \delta_0, (h_\crys^{(l)}-z)^{-1}\delta_{0}\rangle=\kappa^2\langle\chi_{l},(h_l-z)^{-1}\chi_{l}\rangle,\]
\[\kappa_\cS^2\langle\delta_{1},(h_\crys^{(r)}-z)^{-1}\delta_{1}\rangle=\kappa^2\langle\chi_{r},(h_r-z)^{-1}\chi_{r}\rangle.
\]
These relations imply 
\[\kappa_\cS^2 \nu_{\crys}^{(l/r)}=\kappa^2\nu_{l/r},
\]
where $\nu^{(l/r)}_\crys$ is the spectral measure of $h_\crys^{(l/r)}$ associated to $\delta_0/\delta_1$. Since $\nu_{\crys}^{(l/r)}$ and 
$\nu_{l/r}$ are probability measures, we conclude that $\kappa_\cS^2=\kappa^2$ and $\nu_{\crys}^{(l/r)}=\nu_{l/r}$. 
Thus, the transport between the reservoirs is reflectionless iff $\kappa^2=\kappa_\cS^2$ and 
$h_{r/l}$
is unitarily equivalent to $h_\crys^{(r/l)}$. In other words, all suprema in the previous theorem are achieved iff   the EBBM is unitarily equivalent to a crystalline EBMM up 
to (for the transport purposes) irrelevant choice of the sign of $\kappa$.

{\bf Remark 2.} Part~(2) holds whenever $E\Delta_{l/r}(E)$ has a definite sign 
on $\sp(h_\crys)$ and, similarly, Part~(3) holds whenever $\Delta_{l/r}(E)$ has a
definite sign on $\sp(h_\crys)$. In the lack of a definite sign one cannot expect
a variational characterization of Thouless currents in terms of the crystalline
Landauer-B\"uttiker currents. Note that
$$
\Delta_{l/r}(E)=\frac{\sinh((\zeta_{r/l}(E)-\zeta_{l/r}(E))/2)}
{2\cosh(\zeta_l(E)/2)\cosh(\zeta_r(E)/2)},
$$
so that the sign of $\Delta_{l/r}(E)$ is the same as the sign of 
\[
\zeta_{r/l}(E)-\zeta_{l/r}(E)
=(\beta_{r/l}-\beta_{l/r})E-(\beta_{r/l}\mu_{l/r}-\beta_{l/r}\mu_{l/r}).
\]
In the non-trivial case $\beta_l\not=\beta_r$, $\Delta_{l/r}(E)$ changes the 
sign at precisely one point 
\[
E_c=\frac{\beta_{r/l}\mu_{l/r}-\beta_{l/r}\mu_{l/r}}{\beta_{r/l}-\beta_{l/r}}.
\]
If additional information about $\sp(h_\crys)$ is available, the above fact can
be used to obtain further relations between the crystalline Landauer-B\"uttiker 
currents~\eqref{landauer-buttiker-crys} and Thouless currents~\eqref{thoulesscurrents}.

{\bf Remark 3.} If $\beta_l=\beta_r=\infty$, $\mu_l<\mu_r$, and $I=[\mu_l,\mu_r]$, 
then 
\begin{align*}
\langle\Phi_r\rangle_\Th&=\frac{1}{2\pi}\int_{\sp(h_\crys)\cap I}E\d E,\\[3mm]
\langle\cI_r\rangle_\Th&=\frac{1}{2\pi}|\sp(h_\crys)\cap I|.
\end{align*}
Thus, the Thouless formula~\eqref{gTh} indeed describes the maximal conductance at
zero temperature for the given potential interval $I$:
$$
g_\Th(I)=\frac{\langle\cI_r\rangle_\Th}{\mu_r-\mu_l}
=\sup\frac{\langle\cI_r\rangle_\infty}{\mu_r-\mu_l}.
$$

{\bf Remark 4.} Since the crystalline Landauer-B\"uttiker 
formulas~\eqref{landauer-buttiker-crys} are derived 
from the first principles of quantum statistical mechanics, 
Theorem~\ref{final-thouless} can be considered a rigorous quantum statistical 
derivation of the Thouless energy formula. This derivation also 
identifies the heuristic notion of "optimal feeding" of electrons with 
reflectionless transport between the reservoirs.

\
\section{Proof of Theorem  \ref{prop:crys-tra-coef-bound}}
\label{sec-proof}

\subsection{Sample transmittance and Green matrix}

We first connect the transmittance \eqref{tra-coef} to the sample's Green function.
Recall that $F_{l/r}$ are defined in 
\eqref{reserv-funct} and denote by $F(E)$ the $2\times 2$ diagonal matrix with entries $F_l(E)$ and $F_r(E)$. We also introduce the $2\times 2$ Green matrices $G_\cS^{(N)}(z)$ and $G^{(N)}(z)$ with entries
\begin{align*}
G_{\cS,ab}^{(N)}(z)&=\langle\psi_a, (h_{\cS}^{(N)}- z)^{-1}\psi_b\rangle,\\[3mm]
G_{ab}^{(N)}(z)&=\langle\psi_a, (h^{(N)}- z)^{-1}\psi_b\rangle,
\end{align*}
where $a,b\in \{l,r\}$ (recall that $\psi_l=\delta_1$ and $\psi_r=\delta_{NL}$). As usual, we write 
$G_{ab}^{(N)}(E)=G_{ab}^{(N)}(E+\i 0)$.

The full Green matrix $G^{(N)}$ and the sample Green matrix $G_\cS^{(N)}$ are related by (see Lemma 2.1 in \cite{BJP})
$$
G_\cS^{(N)}(E)=(I-\kappa^2 G_\cS^{(N)}(E)F(E))G^{(N)}(E),
$$
from which we deduce
\beq
G_{lr}^{(N)}(E)
=\frac{G_{\cS,lr}^{(N)}(E)}{\det(I-\kappa^2 G_\cS^{(N)}(E)F(E))}.
\label{greenfull-small}
\eeq
Combined with \eqref{tra-coef} and \eqref{landauer-buttiker} this allows us  the expression of 
the transmittance of the sample and hence the steady currents in terms
of the Green matrix $G_\cS^{(N)}$.

\subsection{Green and transfer matrix}

Our next step  is to relate the 
sample's Green matrix to the transfer matrix of the periodic Jacobi matrix $h_\crys$.

Following \cite{Si} the transfer matrix at energy $E$ is defined by
$$
T_n(E) =  A_n(E)\cdots A_1(E),
$$
where 
\[
A_x(E)=J_x^{-1}\left[\begin{array}{cc}
E-\lambda_x & -1 \\
J_x^2 & 0
\end{array}\right].
\]
Note that $\det A_x(E)=1$ for any $x$ and hence $\det T_n(E)=1$ as well.
A function $u$ satisfies the finite difference equation $h_\crys u =Eu$ if and only if 
for any $x$ one has 
$$
A_x\left[\begin{array}{c} u(x) \\ J_{x-1}u(x-1) \end{array} \right]=\left[\begin{array}{c} u(x+1) \\ J_x u(x) \end{array} \right].
$$

\begin{lemma}\label{lem:greentransfer} For any $x,y,u,v\in\cc$ one has
$$
G_\cS^{(N)}(E)\left[\begin{array}{c}x\\y\end{array}\right]
=\left[\begin{array}{c}u\\v\end{array}\right]
\Longleftrightarrow
T_{NL}(E)\left[\begin{array}{c}u\\-x\end{array}\right]
=\left[\begin{array}{c}-\kappa_\cS^{-1}y\\\kappa_\cS v\end{array}\right].
$$
In other words, the matrix $P: (x,y,u,v) \mapsto (u,-x,-\kappa_\cS^{-1}y,\kappa_\cS v)$ maps the graph of $G_\cS^{(N)}(E)$ to that of $T_{NL}(E)$.
\end{lemma}

{\bf Remark.} This lemma is a slight generalization of Lemma 2.2 in \cite{BJP}. We include the  proof for the reader convenience.

\proof Fix $N$ and $E\in\rr\setminus\sp(h_{\cS}^{(N)})$. For $f\in\ell^2(Z_{NL})$, 
the function 
\[
u(x)=\langle\delta_x,(h_{\cS}^{(N)}-E)^{-1}f
\rangle\]
satisfies the finite difference equation
\begin{equation}\label{fde}
(h_\crys-E) u=f,
\end{equation}
with boundary conditions $u(0)=u(NL+1)=0$. Using the transfer matrix
$$
T(x,y)=A_xA_{x-1}\cdots A_{y+1},
$$
the solution of the initial value problem for \eqref{fde} can be written as
$$
\left[\begin{array}{c}u(x+1)\\J_x u(x)\end{array}\right]
=T(x,0)\left[\begin{array}{c}u(1)\\J_0u(0)\end{array}\right]
-\sum_{y=1}^x T(x,y-1)\left[\begin{array}{c}0\\f(y)\end{array}\right].
$$
Setting $x=NL$ and taking the boundary conditions and $J_{NL}=J_0=\kappa_\cS$ into 
account yields
$$
\left[\begin{array}{c}0\\\kappa_\cS u(NL)\end{array}\right]
=T_{NL}(E)\left[\begin{array}{c}u(1)\\0\end{array}\right]
-\sum_{y=1}^{NL}T(NL,y-1)\left[\begin{array}{c}0\\f(y)\end{array}\right],
$$
which is an equation for the unknown $u(1)$ and $u(NL)$. Setting $f=\delta_1$ and 
$f=\delta_{NL}$, we obtain the following equations for the entries of the matrix 
$G_\cS^{(N)}(E)$:
\begin{align*}
T_{NL}(E)\left[\begin{array}{c} G_{\cS,ll}^{(N)}(E) \\ -1 \end{array}\right]
&= \left[\begin{array}{c} 0 \\ \kappa_\cS G_{\cS,rl}^{(N)}(E)\end{array}\right],\\[3mm]
T_{NL}(E)\left[\begin{array}{c} G_{\cS,lr}^{(N)}(E) \\ 0 \end{array}\right]
&= \left[\begin{array}{c} -\kappa_\cS^{-1} \\ \kappa_\cS G_{\cS,rr}^{(N)}(E)\end{array}\right].
\end{align*}
Thus, the two linearly independent vectors 
\[
\left\{[G_{\cS,ll}^{(N)}(E),-1,0,\kappa_\cS G_{\cS,rl}^{(N)}(E)]^T, [G_{\cS,lr}^{(N)}(E),0,-\kappa_\cS^{-1},\kappa_\cS G_{\cS,rr}^{(N)}(E)]^T\right\},
\]
span the graph of $T_{NL}(E)$. One easily checks that
they are the images by the matrix $P$ of the two vectors $[1,0,G_{\cS,ll}^{(N)}(E),G_{\cS,rl}^{(N)}(E)]^T$ and $[0,1,G_{\cS,lr}^{(N)}(E),G_{\cS,rr}^{(N)}(E)]^T$ which span the graph of $G_\cS^{(N)}(E)$.
\hfill \qed


\subsection{Transfer matrix eigenvalues and eigenfunctions}

Since the Jacobi matrix $h_\crys$ is periodic, one has $T_{NL}(E)=T_L(E)^N$ for 
any $N$. The 
eigenvalues and eigenvectors of the one-period transfer matrix $T_L(E)$ will thus 
play an important role. Since $\det T_L(E)=1$, we may write its eigenvalues as 
$\alpha(E)$ and $\alpha(E)^{-1}$. It is a standard result (see, e.g., \cite{RS4,Si})
that these eigenvalues are either complex conjugated (and hence of modulus 1)
or real. The first case occurs iff $|\tr\,T_L(E)|\le2$ which is in turn 
equivalent to $E\in\sp(h_\crys)$. We denote by 
\[
\Psi_{\pm}(E)=\left[\begin{array}{c} 
\phi_{\pm}(E)\\
\kappa_\cS\psi_{\pm}(E)
\end{array}\right]
\]
an eigenvector of $T_L(E)$ associated  to its eigenvalue $\alpha(E)^{\pm1}$ with the following
conventions:
\begin{description}
\item[(N1)\label{N1}] When $E\not\in\sp(h_\crys)$ we chose $|\alpha(E)|^{-1}<1<|\alpha(E)|$ and
real eigenvectors.
\item[(N2)\label{N2}] When $E\in\sp(h_\crys)$ we chose $\Psi_-(E)=\overline{\Psi_+(E)}$.
Since the two eigenvectors are linearly independent, they can be normalized
by $\phi_\pm(E)=1$, which further implies $\Im\psi_\pm(E)\not=0$.
We then select $\alpha(E)$ such that 
$\Im\psi_+(E)>0$.
\end{description}
Using Lemma \ref{lem:greentransfer} and $T_{NL}(E)=T_L(E)^N$, 
we get
$$
G_{\cS}^{(N)}(E)=\frac{-\kappa_\cS^{-1}}{D_N(E)}\left[
\begin{array}{cc}
\phi_+(E)\phi_-(E)(\alpha(E)^N-\alpha(E)^{-N})&\phi_+(E)\psi_-(E)-\phi_-(E)\psi_+(E)\\[6pt]
\phi_+(E)\psi_-(E)-\phi_-(E)\psi_+(E)&\psi_+(E)\psi_-(E)(\alpha(E)^N-\alpha(E)^{-N})
\end{array}
\right]
$$
where 
\[
D_N(E)=\alpha(E)^N\phi_+(E)\psi_-(E)-\alpha(E)^{-N}\phi_-(E)\psi_+(E).
\]
An elementary calculation yields
$$
\det(I-\kappa^2G_\cS^{(N)}(E)F(E))
=\frac{\alpha(E)^N\widetilde\phi_+(E)\widetilde\psi_-(E)
-\alpha(E)^{-N}\widetilde\phi_-(E)\widetilde\psi_+(E)}
{\alpha(E)^N\phi_+(E)\psi_-(E)-\alpha(E)^{-N}\phi_-(E)\psi_+(E)},
$$
where
\beq
\widetilde\psi_\pm(E)=\psi_\pm(E)+\eta^2\kappa_\cS\phi_\pm(E) F_l(E),\qquad
\widetilde\phi_\pm(E)=\phi_\pm(E)+\eta^2\kappa_\cS\psi_\pm(E) F_r(E),
\label{phitilde}
\eeq
and we have set
$$
\eta=\frac{\kappa}{\kappa_\cS}.
$$
Inserting the last relations into \eqref{greenfull-small} leads to the following
expression for the off-diagonal element of the full Green matrix
\beq
G_{lr}^{(N)}(E)
=-\kappa_\cS^{-1}\frac{\phi_+(E)\psi_-(E)-\phi_-(E)\psi_+(E)}
{\alpha(E)^N\widetilde\phi_+(E)\widetilde\psi_-(E)
-\alpha(E)^{-N}\widetilde\phi_-(E)\widetilde\psi_+(E)}.
\label{eq:green-eigendata}
\eeq
\subsection{The large $N$ limit}

We proceed to evaluate the large $N$ limit of $\cT_N(E)$ in distributional sense.

We write the left hand side of~\eqref{distribconverge} as 
$\cT_{N,1}(f)+\cT_{N,2}(f)$ with
\begin{align*}
\cT_{N,1}(f)&=\int_{\sp(h_\crys)\cap\Sigma_l\cap\Sigma_r}\cT_N(E)f(E)\,\d E,\\[4mm]
\cT_{N,2}(f)&=\int_{\rr\setminus(\sp(h_\crys)\cap\Sigma_l\cap\Sigma_r)}
\cT_N(E)f(E)\,\d E.
\end{align*}

To deal with $\cT_{N,2}(f)$ we prove
\begin{lemma}\label{OffSpec}
For almost all $E\in\rr\setminus(\sp(h_\crys)\cap\Sigma_l\cap\Sigma_r)$ one
has
$$
\lim_{N\to\infty}\cT_N(E)=0.
$$
\end{lemma}
\proof
Combining~\eqref{tra-coef} and~\eqref{eq:green-eigendata}, we get
$$
\cT_N(E)=4\kappa^2\eta^2\left|
\frac{\phi_+(E)\psi_-(E)-\phi_-(E)\psi_+(E)}
{\alpha(E)^N\widetilde\phi_+(E)\widetilde\psi_-(E)
-\alpha(E)^{-N}\widetilde\phi_-(E)\widetilde\psi_+(E)}\right|^2\Im F_l(E)\,\Im F_r(E).
$$
It follows from~\eqref{non-tri} that $\cT_N(E)=0$ for almost 
all $E\in\rr\setminus(\Sigma_l\cap\Sigma_r)$ and all $N$. Thus, it suffices to consider
$E\in(\Sigma_l\cap\Sigma_r)\setminus\sp(h_\crys)$. 

For such $E$, the eigenvectors
$\Psi_\pm(E)$ are real and $|\alpha(E)|>1$ by Condition~\ref{N1}. 
Moreover, $\Im F_{l/r}(E)>0$ by~\eqref{ImF}. We claim that this implies
$\widetilde\phi_+(E)\widetilde\psi_-(E)\not=0$ from which the result clearly follows. 
To prove this claim, we argue by contradiction. Assume that $\widetilde\phi_+(E)=0$ (the case 
$\widetilde\psi_-(E)=0$ is similar). It follows from~\eqref{phitilde} that 
$$
\Im\widetilde\phi_+(E)=\eta^2\kappa_\cS\psi_+(E)\Im F_r(E)=0,
$$
and hence $\psi_+(E)=0$. 
Using~\eqref{phitilde} again we get $\phi_+(E)=\widetilde\phi_+(E)=0$. We conclude that
$\Psi_+(E)=0$, a contradiction.\hfill\qed

Since $0\leq \cT_N(E)\leq 1$ and $f\in L^1(\rr)$, Lemma~\ref{OffSpec} and
the dominated convergence theorem yield that $\cT_{N,2}(f)\to0$ as $N\to\infty$.
It thus remains to analyze $\cT_{N,1}(f)$. To this end, we relate the eigenvectors
of $T_L$ to the Weyl $m$-functions~\eqref{Mfunct}.

\medskip
\begin{lemma}\label{lem:eigenfunction-mfunction}
For any $E\in\sp(h_\crys)$, we have
$$
\psi_+(E)=-\frac{1}{\kappa_\cS m_r(E)}, \qquad \psi_-(E)=-\kappa_\cS m_l(E).
$$
\end{lemma}

\proof Let us write the one-period transfer matrix as
\beq
T_L(E)=\left[\begin{array}{cc} a(E) & b(E) \\ c(E) & d(E) \end{array}\right].
\label{TLmatrix}
\eeq
Since $\det\,T_L(E)=1$, we can write the discriminant of the quadratic equation
\beq
c(E)z^2+(a(E)-d(E))z-b(E)=0,
\label{eq:quadratic}
\eeq
as $(\tr\,T_L(E))^2-4$, which is negative for $E\in\sp(h_\crys)$. Thus, 
\eqref{eq:quadratic}
has two complex conjugate solutions. The solution with positive imaginary part is
$m_r(E)$ and the other one is $1/\kappa_\cS^2m_l(E)$. We refer the reader 
to Section~5.2 in~\cite{Si} for a proof of these facts.

Using the representation~\eqref{TLmatrix} and taking Condition~\ref{N2} into account, 
one easily shows that the eigenvalues and eigenvectors of $T_L(E)$ are determined by
$\alpha(E)=\e^{\i\theta(E)}$ and
\beq
\kappa_\cS\psi_\pm(E)=\frac1{b(E)}\left(\e^{\pm\i\theta(E)}-a(E)\right)
=\frac1{b(E)}\left(\frac{d(E)-a(E)}{2}\pm\i\sin\theta(E)\right),
\label{psiform}
\eeq
where $\theta(E)\in]-\pi,\pi[$ is defined by
\beq
\cos\theta(E)=\frac12\tr\,T_L(E),\qquad\sign(\theta(E))=\sign(b(E)).
\label{thetadef}
\eeq
Note that the product of the two solutions of~\eqref{eq:quadratic} is
\beq
-\frac{b(E)}{c(E)}=
\frac{m_r(E)}{\kappa_\cS^2m_l(E)}=|m_r(E)|^2>0,
\label{prodsol}
\eeq
so that $b(E)$ and $c(E)$ have opposite signs. It follows that the solutions of~\eqref{eq:quadratic} are
\begin{align*}
m_r(E)
&=\frac1{c(E)}\left(\frac{d(E)-a(E)}{2}-\i\sin\theta(E)\right),\\[4mm]
\frac1{\kappa_\cS^2m_l(E)}
&=\frac1{c(E)}\left(\frac{d(E)-a(E)}{2}+\i\sin\theta(E)\right).
\end{align*}
Comparing these relations with~\eqref{psiform} and using~\eqref{prodsol} yield the 
result.\hfill\qed

To formulate our next result, let the real functions $r(E)\ge0$ and
$\vartheta(E)$ be defined by the following polar decomposition
\beq
\frac{m_l(E)-\eta^2F_l(E)}{\overline{m}_l(E)-\eta^2F_l(E)}\,
\frac{m_r(E)-\eta^2F_r(E)}{\overline{m}_r(E)-\eta^2F_r(E)}\,
\frac{\overline{m}_r(E)}{m_r(E)}
=r(E)\e^{\i\vartheta(E)}.
\label{rdef}
\eeq

\medskip
\begin{lemma}\label{lem:periodic-limit} Let 
$I\subset\sp(h_\crys)\cap\Sigma_l\cap\Sigma_r$. If there exists $\delta<1$ such that 
$r(E)\le\delta$ for almost all $E\in I$, then
$$
\lim_{N\to\infty}\int_I\cT_N(E)f(E)\,\d E=\int_I\cT_\infty(E)f(E)\,\d E,
$$
for any $f\in L^1(\rr)$.
\end{lemma}

\proof Combining Lemma~\ref{lem:eigenfunction-mfunction} with~\eqref{tra-coef} 
and~\eqref{eq:green-eigendata} we can write the transmittance of the $N$-fold
repeated sample as
\beq
\cT_N(E)=\frac{16\eta^4\Im m_r(E)\,\Im m_l(E)}
{|\overline{m}_r(E)-\eta^2F_r(E)|^2|\overline{m}_l(E)-\eta^2F_l(E)|^2}
\frac{\Im\,F_l(E)\,\Im\,F_r(E)}{|1-r(E)\e^{\i(2N\theta(E)+\vartheta(E))}|^2},
\label{TNform}
\eeq
where $\theta(E)$ is defined by~\eqref{thetadef}.
Expanding the right hand side of~\eqref{TNform} in powers of $r(E)$, one obtains
$$
\cT_N(E)=\cT_\infty(E)
\sum_{k\in\zz}r(E)^{|k|}\e^{\i k(2N\theta(E)+\vartheta(E))},
$$
where $\cT_\infty(E)$ is given by~\eqref{def:crys-tra-coef}. Since $r(E)\le\delta<1$
on $I$, this expansion is uniformly convergent for $E\in I$, and we have
$$
\lim_{N\to\infty}\int_I\cT_N(E)f(E)\,\d E 
=\lim_{N\to\infty}\sum_{k\in\zz}\int_I\cT_\infty(E)r(E)^{|k|}
\e^{\i k(2N\theta(E)+\vartheta(E))}f(E)\,\d E.
$$
Since the  function $\theta(E)$ is strictly monotone in each band of $h_\crys$ (see, e.g., Sections~5.3-5.4 in~\cite{Si}), the Riemann-Lebesgue lemma yields the result.
\hfill\qed

For $k>0$, set
$$
I_k=\{E\in\sp(h_\crys)\cap\Sigma_l\cap\Sigma_r\,|\,
\Im F_{l/r}(E)\geq1/k \ \mbox{ and } \ |\tr\,T_L(E)|<2-1/k\},
$$
and $I'_k=(\sp(h_\crys)\cap\Sigma_l\cap\Sigma_r)\setminus I_k$. Obviously, 
\beq
\lim_{k\to\infty}|I'_k|=0.
\label{residual}
\eeq
It follows from the proof of Lemma~\ref{lem:eigenfunction-mfunction} that there
is $\epsilon_k>0$ such that $\Im m_{l/r}(E)\geq\epsilon_k$ for
almost every $E\in I_k$. One easily concludes that there exists
$\delta_k<1$ such that
$$
\left|\frac{m_{l/r}(E)-\eta^2F_{l/r}(E)}{\overline{m}_{l/r}(E)-\eta^2F_{l/r}(E)}\right|
\le\delta_k,
$$
holds for almost every $E\in I_k$. Thus, for such $E$,
$$
r(E)=\left|\frac{m_l(E)-\eta^2F_l(E)}{\overline{m}_l(E)-\eta^2F_l(E)}\right|\,
\left|\frac{m_r(E)-\eta^2F_r(E)}{\overline{m}_r(E)-\eta^2F_r(E)}\right|\le\delta_k^2<1.
$$

Writing
\begin{align*}
\left|\cT_{N,1}(f)-\int\cT_\infty(E)f(E)\,\d E\right|&\le
\left|\int_{I_k}(\cT_N(E)-\cT_\infty(E))f(E)\,\d E\right|\\[3mm]
&+\int_{I'_k}|(\cT_N(E)-\cT_\infty(E))f(E)|\,\d E,
\end{align*}
and applying Lemma~\ref{lem:periodic-limit} to $I_k$, we get
$$
\limsup_{N\to\infty}\left|\cT_{N,1}(f)-\int\cT_\infty(E)f(E)\,\d E\right|
\le 2\int_{I'_k}|f(E)|\,\d E.
$$
This estimate and~\eqref{residual} yield 
\[
\lim_{N\to\infty}\cT_{N,1}(f)=\int\cT_\infty(E)f(E)\,\d E
\]
and Theorem~\ref{prop:crys-tra-coef-bound} follows.


\end{document}